\documentclass[12pt]{article}
\usepackage[russian, english]{babel}
\usepackage{srcltx}
\usepackage{amssymb}
\usepackage{amsfonts, euscript}
\usepackage{amsmath}
\usepackage{amsthm}
\usepackage{enumerate}
\usepackage[cp1251]{inputenc}
\usepackage[T2A]{fontenc}
\usepackage{multicol}

\begin{document}

\selectlanguage{russian}

\vspace*{5mm} \noindent {\sl УДК}: 517.54

\begin{center}
ПОЛНАЯ ТОПОЛОГИЧЕСКАЯ КЛАССИФИКАЦИЯ ПРОСТРАНСТВ БЭРОВСКИХ 
ФУНКЦИЙ НА ОРДИНАЛАХ\footnote{Работа выполнена всеми авторами при финансовой поддержке Российского фонда 	фундаментальных исследований (код проекта 17–51–18051)}
\end{center}

\begin{center}
\textbf{Л.В.~Гензе, С.П.Гулько, Т.Е.~Хмылева}
\end{center}

\noindent\textbf{Аннотация.} В работе рассматриваются пространства
$B_p[1,\alpha]$ 
всех бэровских 
функций 
$x\colon [1,\alpha]\to\mathbb{R}$,
определённых на отрезках ординалов
$[1,\alpha]$
и наделённых топологией поточечной сходимости. Даётся полная топологическая классификация 
этих пространств

\vspace*{5mm} \noindent\textbf{Ключевые слова:} функция 1-го класса Бэра, пространство бэровских функций, топология поточечной сходимости, гомеоморфизм, отрезок ординалов, порядковая топология, вещественная компактность.

\begin{center}
\textbf{\S1. Введение}
\end{center}

Полная линейная топологическая классификация банаховых прост\-ранств $C[1,\alpha]$ всех непрерывных функций на компактных отрезках ординалов была проведена в работах \cite{BesPel,Semadeni,GulOs,Kisl}. Линейную топологическую классификацию этих пространств, но уже  в топологии поточечной сходимости (т.е. пространств   $C_p[1,\alpha]$) см. в \cite{Gulko, BaarsDeGroot}. Аналогичная линейная топологическая классификация для пространств бэровских функций $B_p[1,\alpha]$ также наделенных топологией поточечной сходимости была проведена в работе \cite{GGK}.
В данной работе дана полная  топологическая классификация пространств $B_p[1,\alpha]$ (См. теорему 2.2 ниже). Оказалось, что линейная топологическая классификация совпадает с топологической классификацией.

Наша терминология в основном следует \cite{Engel}.
Отрезки ординалов 
$[1,\alpha]$
наделяются порядковой топологией $\Im$.
Если
$\alpha$ ---
произвольный ординал, а
$\lambda$ ---
начальный ординал (т.е. кардинал), не превосходящий 
$\alpha$,
то положим
$$A_{\lambda,\alpha}=\{t\in [1,\alpha] : \chi (t)=|\lambda|\},$$
где
$\chi (t)$ ---
характер точки 
$t\in [1,\alpha]$.
В частности, 
$A_{\omega,\alpha}$ ---
это множество всех тех точек  
$t\in [1,\alpha]$,
для которых 
$\chi (t)=\aleph_0$.

Символом
$\Im_\omega$
будем обозначать $\aleph_0$~--~модификацию топологии $\Im$, т.е. такую топологию, в которой открытыми объявляются все $G_\delta$~- множества (т.е. те множества, которые можно представить в виде пересечения счётного семейства элементов $\Im$).
Отрезок ординалов 
$[1,\alpha]$,
наделённый топологией
$\Im_\omega$,
будем обозначать 
$[1,\alpha]_\omega$.

Отметим следующие свойства пространства
$[1,\alpha]_\omega$:

($a_\omega$) Если
$A$ --- 
счётное подмножество в
$[1,\alpha]_\omega$,
то
$A$
замкнуто и дискретно. 

($b_\omega$) Для любого ординала 
$\alpha$ 
пространство 
$[1,\alpha]_\omega$
линделёфово
и, следовательно, нормально.

Действительно, это очевидно, если 
$\alpha\leqslant\omega_1$,
где
$\omega_1$ --- 
первый несчётный ординал.
Для 
$\alpha >\omega_1$
линделёфовость пространства 
$[1,\alpha]_\omega$
несложно доказать методом трансфинитной индукции. 

Пусть 
$\alpha$ ---
предельный ординал. Наименьший порядковый тип множеств 
$A\subset [1,\alpha]$,
конфинальных в
$[1,\alpha)$,
будем называть \textit{конфинальностью} ординала 
$\alpha$
и обозначать 
$\mathrm{cf}(\alpha)$.

Нетрудно видеть, что $|\mathrm{cf}(\alpha)|=\chi(\alpha)$ для предельного ординала $\alpha$.
Начальный ординал 
$\alpha$
называется \textit{регулярным}, если
$\mathrm{cf}(\alpha) = \alpha$. В противном случае начальный ординал называется \textit{сингулярным}.
 
\textbf{Определение 1.1.} 
Функцию 
$x\colon [1,\alpha]\to\mathbb{R}$
будем называть \textit{функцией 1-го класса Бэра}, если существует последовательность непрерывных функций
$x_n\colon [1,\alpha]\to\mathbb{R}$,
поточечно сходящаяся к функции
$x$.
Множество всех функций 1-го класса Бэра будем обозначать 
$B^1[1,\alpha]$.
 Если $\gamma>1$ есть счетный ординал, то обозначим $B^{\gamma} [1,\alpha]$ совокупность всех функций, представимых в виде поточечного предела последовательности функций $f_n:[1,\alpha]\to\mathbb R$, где $f_n\in B^{\beta_n}[1,\alpha]$ и $\beta_n<\gamma$ --- счетные ординалы. Это же множество, снабжённое топологией поточечной сходимости, будем обозначать
$B_p^\gamma [1,\alpha]$.

\textbf{Предложение 1.2.} Функция 
$x\colon [1,\alpha]\to\mathbb{R}$ 
является функцией 1-го класса Бэра тогда и только тогда, когда она непрерывна во всех точках
$t\in [1,\alpha]\setminus A_{\omega,\alpha}$.

\begin{proof}
Пусть 
$x\in B_p [1,\alpha]$. 
Тогда существует последовательность функций 
$x_n \in C[1,\alpha]$, 
сходящаяся к 
$x$ 
в каждой точке. Зафиксируем 
$\beta\in [1,\alpha]$ 
со свойством $\mathrm{cf}(\beta)>\omega$. 
Тогда для каждого натурального 
$n$ 
найдётся такой ординал 
$\beta_n < \beta$, 
что 
$x_n(\gamma) = x_n(\beta)$ 
при всех 
$\gamma\in (\beta_n,\beta]$ 
\cite[с. 206]{Engel} (другими словами, каждое отображение $x_n$ становится постоянным в некоторой окрестности точки $\beta$). 
Пусть 
$\beta_0 = \sup\{\beta_n : n\in \mathbb{N}\}$. 
Так как 
$\mathrm{cf}(\beta)>\omega$, 
то 
$\beta_0 < \beta$ 
и при любом 
$\gamma\in (\beta_0,\beta]$ 
выполнено
$$
x(\gamma)=\lim_{n\to \infty}x_n(\gamma)=\lim_{n\to \infty}x_n(\beta)=x(\beta).
$$ 
Следовательно, отображение 
$x$ 
непрерывно в точке 
$\beta$.

Обратно, пусть отображение 
$x\colon [1,\alpha]\to \mathbb{R}$ 
непрерывно во всех точках несчётной конфинальности. Докажем, что найдётся последовательность непрерывных отображений 
$\{x_n : n\in \mathbb{N}\}$, 
поточечно сходящихся к 
$x$. 
Доказательство проведём по трансфинитной индукции. Ясно, что если 
$\alpha$ --- 
конечный ординал, то утверждение теоремы верно.

Предположим, что для всех ординалов, меньших 
$\alpha$, 
утверждение уже доказано.

Случай 1:
 $\alpha$ --- 
бесконечный непредельный ординал. По индуктивному предположению существует такая последовательность непрерывных функций 
$\{y_n : n\in \mathbb{N}\}$, 
что 
$\displaystyle\lim_{n\to \infty} y_n(\gamma)=x(\gamma)$ 
для всех 
$\gamma\in [1,\alpha-1]$. 
Продолжим функции 
$y_n$ 
на отрезок 
$[1,\alpha]$, 
полагая
$$
x_n(\gamma) = \left\{
\begin{array}{rl}
y_n(\gamma), & \gamma\in [1,\alpha-1]; \\
x(\alpha), & \gamma = \alpha. \\
\end{array} \right. 
$$
Очевидно, что 
$\{x_n : n\in \mathbb{N}\}$ --- 
требуемая последовательность непрерывных функций.

Случай 2: 
$\alpha$ --- 
предельный ординал и 
$\mathrm{cf}(\alpha)=\omega$. 
Тогда существует возрастающая последовательность ординалов 
$\{\alpha_n : n\in \mathbb{N}\}$, 
такая, что 
$\lim_{n\to \infty}\alpha_n=\alpha$. 
Рассмотрим следующее разбиение отрезка 
$[1,\alpha]$:
$$
[1,\alpha] = [1,\alpha_1]\cup(\alpha_1,\alpha_2]\cup \ldots \cup(\alpha_{k-1},\alpha_k]\cup \ldots \cup\{\alpha\}.
$$
На отрезке 
$[1,\alpha_1]$ 
и на каждом из полуинтервалов 
$(\alpha_{k-1},\alpha_k]$ 
существует последовательность непрерывных функций 
$\{y_n^k : n\in \mathbb{N}\}$, 
поточечно сходящаяся к 
$x$. 
Положим
$$
x_n(\gamma) = \left\{
\begin{array}{rl}
y_n^1(\gamma), & \gamma\in [1,\alpha_1]; \\
y_n^2(\gamma), & \gamma\in (\alpha_1,\alpha_2]; \\
\vdots \\
y_n^n(\gamma), & \gamma\in (\alpha_{n-1},\alpha_n]; \\
x(\alpha), & \gamma\in (\alpha_n,\alpha]. \\
\end{array} \right. 
$$
Ясно, что все функции 
$x_n$ 
непрерывны на отрезке 
$[1,\alpha]$ 
и поточечно сходятся на этом отрезке к 
$x$.

Случай 3: 
$\alpha$ --- 
предельный ординал и 
$\mathrm{cf}(\alpha)>\omega$. 
Тогда существует такой ординал 
$\gamma_0<\alpha$, 
что 
$x(\gamma)=x(\alpha)$ 
при всех 
$\gamma\in (\gamma_0,\alpha]$. 
По предположению индукции существует последовательность непрерывных функций 
$\{y_n : n\in \mathbb{N}\}$, 
заданных на отрезке $[1,\gamma_0]$
и поточечно сходящихся к 
$x$. 
Продолжим эти функции на отрезок 
$[1,\alpha]$, 
полагая
$$
x_n(\gamma) = \left\{
\begin{array}{rl}
y_n(\gamma), & \gamma\in [1,\gamma_0]; \\
x(\alpha), & \gamma\in (\gamma_0,\alpha]. \\
\end{array} \right. 
$$
Понятно, что 
$\{x_n : n\in \mathbb{N}\}$ 
сходится поточечно к 
$x$
на отрезке
$[1,\alpha]$.
\end{proof}

Из доказанного критерия нетрудно вывести, что поточечный предел последовательности функций первого класса также будет функцией первого класса, следовательно, 
$B_p^\gamma [1,\alpha]=B_p^1 [1,\alpha]$ для каждого счетного ординала $\gamma$.
Это означает, что на отрезке ординалов $[1,\alpha]$ пространство всех бэровских функций совпадает с с пространством функций первого класса.  В дальнейшем мы будем обозначать это пространство через $B_p[1,\alpha]$.
При таком соглашении из предложения 1.2 следует:

\textbf{Следствие 1.3.} 
$B_p[1,\alpha]=C_p([1,\alpha]_\omega)$.

\textbf{Следствие 1.4.} Если
$\alpha < \omega_1$,
то любая функция 
$x\colon [1,\alpha]\to\mathbb{R}$
является функцией первого класса Бэра.

\textbf{Предложение 1.5.} Функция 
$x\colon [1,\alpha]\to\mathbb{R}$
является 
является бэровской
тогда и только тогда, когда множество 
$x^{-1}(U)$ 
имеет тип
$F_\sigma$
для любого множества 
$U$,
открытого в 
$\mathbb{R}$.
\begin{proof}
Пусть 
$x\colon [1,\alpha]\to\mathbb{R}$ 
и множество 
$x^{-1}(U)$ 
имеет тип 
$F_\sigma$ 
в 
$[1,\alpha]$ 
для любого открытого 
$U\subset\mathbb{R}$. 
Это эквивалентно тому, что 
$x^{-1}(H)$ 
имеет тип 
$G_\delta$ 
в 
$[1,\alpha]$ 
для любого замкнутого множества 
$H\subset\mathbb{R}$. 
В силу следствия 1.4 можно считать, что 
$\alpha\geqslant\omega_1$. 
Рассмотрим ординал 
$\beta\leqslant\alpha$ 
такой, что 
$\mathrm{cf}(\beta)>\omega$. 
Тогда 
$x^{-1}\left(x(\beta)\right)$ --- 
это 
$G_\delta$-множество, 
содержащее точку 
$\beta$, 
следовательно, содержащее некоторый полуинтервал 
$(\gamma, \beta]$. 
Таким образом, на множестве 
$(\gamma, \beta]$ 
отображение 
$x$ 
постоянно, а значит, непрерывно в точке 
$\beta$. 
По предложению 1.2 заключаем, что 
$x$ --- 
отображение первого класса Бэра в смысле определения 1.1.  

Обратно, пусть 
$U$ --- 
открытое подмножество в 
$\mathbb{R}$. 
Представим 
$U$ 
в виде 
$U=\bigcup\limits_{k=1}^{\infty} F_k$, 
где все 
$F_k$ 
замкнуты в 
$\mathbb{R}$ 
и 
$F_k \subset \mathrm{Int}\, F_{k+1}$, $k=1,2,\ldots$. 
Рассмотрим последовательность непрерывных функций 
$x_n\colon [1,\alpha]\to \mathbb{R}$, 
поточечно сходящуюся к 
$x$. 
Нетрудно убедиться, что
\begin{equation*}
x^{-1}(U)=\bigcup_{k=1}^{\infty} \bigcap_{n=k}^{\infty} x_n^{-1}(F_k), 
\end{equation*}
откуда следует, что 
$x^{-1}(U)$ есть 
$F_\sigma$-множество в 
$[1,\alpha]$.
\end{proof}

\textbf{Определение 1.6.} (\cite{Arch}) Пусть
$Z$ ---
топологическое пространство. Функция
$x\colon Z\to\mathbb{R}$
называется \textit{строго $\omega$-непрерывной}, если для любого счётного множества 
$D\subset Z$
найдётся такая непрерывная функция
$y\colon Z\to\mathbb{R}$,
что 
$y|_D = x|_D$.

\textbf{Предложение 1.7.} Любая функция 
$x\colon [1,\alpha]_\omega\to\mathbb{R}$
является строго 
$\omega$-непрерывной.
\begin{proof}
Пусть
$D\subset [1,\alpha]_\omega$ ---
счётное подмножество. По свойству ($a_\omega$) функция $x|_D$ непрерывна на замкнутом множестве $D$.  
Кроме того, по свойству ($b_\omega$) и
по теореме Титце--Урысона существует непрерывное продолжение
$y: [1,\alpha]_\omega\to\mathbb{R}$.
\end{proof}

Напомним \cite{Engel}, что топологическое пространство $X$ называется вещественно полным, если оно гомеоморфно замкнутому подпространству некоторого произведения вещественных прямых. Вещественно полное пространство $\nu X$ называется вещественной компактификацией пространства $X$, если существует гомеоморфное вложение $\nu:X\to \nu X$, для которого замыкание ${\overline{\nu (X)}}$ совпадает с $\nu X$ и для любой непрерывной функции $f:X\to \mathbb{R}$ найдется непрерывное продолжение $\tilde f:\nu X\to\mathbb{R}$ .

\textbf{Предложение 1.8.} Вещественная компактификация
$\nu B_p[1,\alpha]$
канонически гомеоморфна тихоновскому произведению
$\mathbb{R}^{[1,\alpha]}$.

\begin{proof}
Поскольку для любого топологического пространства
$X$
пространство
$\nu C_p(X)$
канонически гомеоморфно множеству всех строго 
$\omega$-непрерывных функций на
$X$
(\cite{Tkachuk}, стр.382), то, учитывая следствие 1.3 и предложение 1.7, получаем
$\nu B_p [1,\alpha] = \nu C_p \bigl([1,\alpha]_\omega\bigr) = \mathbb{R}^{[1,\alpha]}$.
\end{proof}

\begin{center}
\textbf{\S2. Доказательство основной теоремы}
\end{center}

В работе \cite{GGK} доказана следующая теорема.

\textbf{Теорема 2.1.} Пусть
$\alpha$ и $\beta$ ---
бесконечные ординалы и 
$\alpha \leqslant \beta$.
Тогда пространства 
$B_p[1,\alpha]$ и $B_p[1,\beta]$,
наделённые топологиями поточечной сходимости, линейно гомеоморфны тогда и только тогда, когда выполняется одно из следующих взаимоисключающих условий:
\begin{enumerate}
	\item $\omega\leqslant\alpha\leqslant\beta <\omega_1$;
	\item $\omega_1\leqslant\alpha\leqslant\beta <\omega_2$;
	\item $\tau\cdot n\leqslant\alpha\leqslant\beta <\tau\cdot (n+1)$, где $\tau\geqslant\omega_2$ --- начальный регулярный ординал и $n < \omega$;
	\item $\tau\cdot\sigma\leqslant\alpha\leqslant\beta <\tau\cdot \sigma^+$, где $\tau\geqslant\omega_2$ --- начальный регулярный ординал, $\sigma$ --- такой начальный ординал, что $\omega\leqslant\sigma < \tau$ и $\sigma^+$ --- наименьший начальный ординал, больший, чем $\sigma$;
	\item $\tau^2\leqslant\alpha\leqslant\beta <\tau^+$, где $\tau\geqslant\omega_2$ --- начальный регулярный ординал и $\tau^+$ --- наименьший начальный ординал, больший, чем $\tau$;
	\item $\tau\leqslant\alpha\leqslant\beta <\tau^+$, где $\tau\geqslant\omega_2$ --- начальный сингулярный ординал.
\end{enumerate} 

В данной работе мы доказываем следующую теорему:

\textbf{Теорема 2.2.} Пусть
$\alpha$ и $\beta$ ---
бесконечные ординалы.
Пространства 
$B_p[1,\alpha]$ и $B_p[1,\beta]$
гомеоморфны тогда и только тогда, когда они линейно гомеоморфны.

Ясно, что если 
$\alpha$ и $\beta$
попадают в один из указанных в теореме~1 промежутков, то 
$B_p[1,\alpha]$ и $B_p[1,\beta]$
гомеоморфны. Если же 
$\alpha$ и $\beta$
попадают в разные промежутки и 
$\alpha <\tau\leqslant\beta$
для некоторого начального ординала 
$\tau$,
то пространства
$B_p[1,\alpha]$ и $B_p[1,\beta]$
не гомеоморфны, так как эти пространства имеют различный вес.

Таким образом, чтобы доказать теорему 2.2, нам достаточно доказать, что пространства
$B_p[1,\tau\cdot\sigma]$ и $B_p[1,\tau\cdot\lambda]$
не гомеоморфны, если 
$\tau\geqslant\omega_2$
 --- начальный регулярный ординал, а  
$\sigma$ и $\lambda$ ---
начальные ординалы, удовлетворяющие условию
$1\leqslant\sigma <\lambda\leqslant\tau$.

\textbf{Лемма 2.3.} Пусть 
$\alpha$ --- 
произвольный ординал, $\lambda$ ---
начальный регулярный ординал, $\omega_1<\lambda \leq \alpha$,
функция
$x\colon [1,\alpha]\to\mathbb{R}$
непрерывна в точках множества 
$A_{\omega_1,\alpha}$ и $t_0\in A_{\lambda,\alpha}$.
Тогда существует такой ординал
$\gamma < t_0$,
что 
$x|_{(\gamma,t_0)}=\mathrm{const}$.

\begin{proof}
Предположим противное. Так как 
$\mathrm{cf}(t_0)=\lambda > \omega_1$,
то для некоторого 
$\varepsilon_0 >0$ 
и для каждого
$\gamma < t_0$
найдутся ординалы
$t_\gamma$
и
$q_\gamma$,
удовлетворяющие неравенству
$\gamma < t_\gamma < q_\gamma < t_0$
и такие, что
$|x(t_\gamma) - x(q_\gamma)|\geqslant \varepsilon_0$.

По трансфинитной индукции для каждого
$\xi\in [1,\omega_1)$
можно выбрать точки
$\gamma_\xi$, $t_{\gamma_\xi}$
и
$q_{\gamma_\xi}$
так, чтобы
$\gamma_\xi < t_{\gamma_\xi} < q_{\gamma_\xi} < \gamma_{\xi +1} < t_0$
и 
$\gamma_\xi = \sup_{\eta < \xi} \gamma_\eta$
для предельного ординала 
$\xi$.
Тогда точка
$\gamma_0 = \sup_{\xi\in [1,\omega_1)}\gamma_\xi = \sup_{\xi\in [1,\omega_1)}q_{\gamma_\xi} = \sup_{\xi\in [1,\omega_1)} t_{\gamma_\xi}$
будет элементом множества 
$A_{\omega_1,\alpha}$,
а функция 
$x$
будет разрывна в точке
$\gamma_0$,
что противоречит условию леммы.
\end{proof}

Введём некоторые обозначения. Для функции
$x\in \mathbb{R}^{[1,\alpha]}$
и начального ординала
$\lambda\leqslant\alpha$
символом 
$G_\lambda (x)$
будем обозначать семейство

\begin{multline*}
G_\lambda (x) = \left\{ \bigcap_{s\in S}V_s : V_s\text{ --- стандартная окрестность } x \text{ в }\mathbb{R}^{[1,\alpha]}\text{ и }  |S| = \lambda \right\}.
\end{multline*}

Элементы семейства
$G_\lambda (x)$
будем называть 
\textit{$\lambda$-окрестностями функции $x$}.

Символом 
$D(x)$
будем обозначать множество всех точек разрыва функции
$x\in \mathbb{R}^{[1,\alpha]}$.

Для регулярного ординала
$\tau > \omega_1$
и начального ординала
$\sigma$ такого, что $\sigma\leqslant\tau$
положим
\begin{multline*}
M_{\tau\sigma} = \left\lbrace x\in \mathbb{R}^{[1,\tau\cdot\sigma]} : x \text{ непрерывна в тех точках }t\in [1,\tau\cdot\sigma],\right. \\
\left. \vphantom{\mathbb{R}^{[1,\tau\cdot\sigma]}}\text{ для которых } \omega_1 \leqslant \mathrm{cf}(t)<\tau\right\rbrace.
\end{multline*}

Ясно, что 
$B_p[1,\tau\cdot\sigma]\subset M_{\tau\sigma}$.

\textbf{Лемма 2.4.} Пусть
$\tau > \omega_1$ --- 
начальный регулярный ординал,
$\sigma\leqslant \tau$ ---
начальный ординал и 
$x\in M_{\tau\sigma}$.
Тогда множество
$D(x)$
не более чем счётно.

\begin{proof}
Предположим, что множество 
$D(x)$
несчётно и пусть
$$
t_0 = \min\{t\in [1,\tau\cdot\sigma] : [1,t]\cap D(x)\text{ несчётно}\}.
$$
Ясно, что
$\mathrm{cf}(t_0) > \omega$. 
Так как 
$x\in M_{\tau\sigma}$,
то
$x$
непрерывна в тех точках
$t$,
для которых
$\mathrm{cf}(t)=\omega_1$. 
По лемме 2.3 найдётся ординал   
$\gamma < t_0$ такой, что
$x|_{(\gamma,t_0)}$ ---
постоянная функция (если
$\mathrm{cf}(t_0) = \omega_1$,
то такой ординал
$\gamma$
существует в силу непрерывности функции
$x$
в точке
$t_0$).
  
Но тогда, с учётом равенства 
$$
[1,t_0]\cap D(x) = ([1,\gamma]\cap D(x))\cup ([\gamma,t_0]\cap D(x)),
$$
множество
$[1,\gamma]\cap D(x)$
несчётно, что противоречит определению ординала
$t_0$.
\end{proof}

\textbf{Лемма 2.5.} Пусть 
$\tau > \omega_1$ --- 
начальный регулярный ординал,
$\sigma$
 --- начальный ординал и
$\sigma\leqslant\tau$. 
Тогда верна формула
\begin{multline*}
M_{\tau\sigma} = \left\lbrace x\in \mathbb{R}^{[1,\tau\cdot\sigma]} : V\cap B_p[1,\tau\cdot\sigma]\ne\varnothing  \text{ для каждого начального}\right. \\
\left.\vphantom{\mathbb{R}^{[1,\tau\cdot\sigma]}} \text{ординала }\lambda < \tau 
\text{ и любой $\lambda$-окрестности $V$ функции } x \right\rbrace.
\end{multline*}
\begin{proof}
Обозначим правую часть равенства 
$L_{\tau\sigma}$ и докажем, что это множество совпадает с $M_{\tau\sigma}$.
Предположим, что 
$x\notin M_{\tau\sigma}$,
то есть 
$x$
разрывна в некоторой точке
$t_0$,
для которой 
$\omega_1 \leqslant \mathrm{cf}(t_0) < \tau$.
Поскольку 
$|\mathrm{cf}(t_0)| = \chi(t_0)$,
то существует такая база
$\{U_j(t_0)\}_{j\in J}$
окрестностей точки
$t_0$,
что 
$|J|<\tau$.
Так как функция 
$x$
разрывна в точке
$t_0$,
то существует такое число
$\varepsilon_0 > 0$,
что для каждого
$j\in J$
найдётся точка
$t_j\in U_j(t_0)$,
для которой
$|x(t_j)-x(t_0)|\geqslant\varepsilon_0$.
Пусть 
$V = \cap_{j\in J, n\in\mathbb{N}}V_{j,n}$,
где
$V_{j,n} = V(x,t_j,t_0,1/n)$ ---
стандартная окрестность функции
$x$
в пространстве 
$\mathbb{R}^{[1,\tau\cdot\sigma]}$.
Если
$y\in V$,
то
$y(t_j)=x(t_j)$
и
$y(t_0)=x(t_0)$.
Следовательно, функция
$y$
разрывна в точке 
$t_0$
и тогда
$y\notin B_p [1,\tau\cdot\sigma]$.
Таким образом,
$V\cap B_p[1,\tau\cdot\sigma] = \varnothing$,
то есть
$x\notin L_{\tau\sigma}$.

Пусть теперь 
$x\in M_{\tau\sigma}$, т.е.  функция $x$ может быть разрывна только в точках множества 
$A_{\tau,\tau\cdot\sigma}$ и в точках, конфинальных $\omega$.

Несложно убедиться в том, что множество 
$A_{\tau,\tau\cdot\sigma}$ имеет вид
$$
A_{\tau,\tau\cdot\sigma}=\{\tau\cdot (\xi+1) : \xi <\sigma\},  \text{ если } \sigma <\tau;
$$
$$
A_{\tau,\tau\cdot\sigma}=\{\tau\cdot (\xi+1) : \xi <\tau\}\cup \{\tau\cdot\tau\}, \text{ если } \sigma =\tau. 
$$

По лемме 2.4 множество
$D(x)$
не более чем счётно и, следовательно, для некоторой последовательности $\{\xi_n\}_{n=1}^\infty\subset [1,\sigma)$
$$
A_{\tau,\tau\cdot\sigma}\cap D(x)=
\{\tau \cdot(\xi_n+1) : n\in\mathbb{N}\},  \text{ если } \sigma <\tau.
$$
Аналогично, 
$$
A_{\tau,\tau\cdot\sigma}\cap D(x)=
\{\tau\cdot(\xi_n+1) : n\in\mathbb{N}\}\cup \{\tau\cdot\tau\},  \text{ если } \sigma =\tau.
$$


Пусть $\lambda <\tau$ и 
$V(x) = \bigcap\{U(x,\eta,1/n): \eta\in S,\,n\in\mathbb{N}\}$ --- $\lambda$-окрестность точки 
$x$.
Тогда 
$|S| <\tau$.

Так как  множество $S$
не конфинально регулярному ординалу
$\tau$,
то для каждого
$n\in\mathbb{N}$
найдётся ординал
$\gamma_n$
такой, что
$\tau\xi_n <\gamma_n < \tau (\xi_n+1)$
и
$(\gamma_n,\tau(\xi_n+1))\cap S = \varnothing$.

В случае
$\sigma = \tau$
также найдётся такой ординал
$\gamma_0 < \tau^2$,
что 
$(\gamma_0,\tau^2)\cap S = \varnothing$
и
$(\gamma_0,\tau^2)\cap\{\tau(\xi_n+1)\}_{n=1}^\infty = \varnothing$.

Рассмотрим функцию
$$
\tilde{x}(t)=\left\lbrace 
\begin{array}{ll}
x(\tau(\xi_n+1)), & \text{ если } t\in (\gamma_n, \tau(\xi_n+1));\\
x(\tau^2), & \text{ если } t\in (\gamma_0, \tau^2);\\ x(t), & \text{ в остальных случаях}.
\end{array}\right. 
$$
Так как 
$\tilde{x}|_S = x|_S$,
то 
$\tilde{x}\in V(x)$.
С другой стороны,
пусть $t\in [1,\tau\cdot\sigma]$  и $\mathrm{cf}(t)\geqslant\omega_1$. Очевидно, что функция $\tilde{x}$ постоянна на каждом из промежутков $(\gamma_n, \tau(\xi_n+1)]$ и $(\gamma_0, \tau^2]$ и, следовательно, непрерывна во всех точках $t\in\cup_{n=1}^\infty(\gamma_n, \tau(\xi_n+1)]\cup (\gamma_0, \tau^2]$. Если $t\notin \overline {\cup_{n=1}^\infty(\gamma_n, \tau(\xi_n+1)]}\cup (\gamma_0, \tau^2]$, то функция $\tilde{x}$
совпадает с функцией
$x$
в некоторой окрестности точки 
$t$ и, значит, непрерывна в $t$. Если же $t\in \overline {\cup_{n=1}^\infty(\gamma_n, \tau(\xi_n+1)]}\backslash {\cup_{n=1}^\infty(\gamma_n, \tau(\xi_n+1)]}$, то $\mathrm{cf}(t)=\omega$. Таким образом 
$\tilde{x}\in B_p [1,\tau\cdot\sigma]$,
то есть
$V(x)\cap B_p [1,\tau\cdot\sigma]\neq \varnothing$
и, следовательно, 
$x\in L_{\tau\sigma}$.
\end{proof}

Доказательство следующей леммы можно найти в (\cite{Engel}, стр.327).

\textbf{Лемма 2.6.} Пусть $X$ и $Y$ --- топологические пространства и  
$\varphi\colon X\to Y$ --- 
гомеоморфизм. Тогда  существует гомеоморфизм
$\tilde{\varphi}\colon \nu X\to \nu Y$ такой,
что
$\tilde{\varphi}(x) = \varphi (x)$
для каждого
$x\in X$. $\Box$

\textbf{Предложение 2.7.} Пусть  
$\tau > \omega_1$ ---
регулярный ординал,
$\sigma, \eta$ ---
такие начальные ординалы, что
$\omega\leq\eta <\sigma\leqslant\tau$.
Тогда пространства
$B_p [1,\tau\cdot\sigma]$
и
$B_p [1,\tau\cdot\eta]$
не гомеоморфны.
\begin{proof}
Предположим, что существует гомеоморфизм
$$\varphi\colon B_p [1,\tau\cdot\sigma]\to B_p [1,\tau\cdot\eta].
$$
Не нарушая общности, можно считать, что $\varphi(0)=0$. По лемме 2.6 существует гомеоморфизм
$$
\tilde{\varphi}\colon \nu (B_p [1,\tau\cdot\sigma])\to \nu(B_p [1,\tau\cdot\eta])
$$
такой, что 
$\tilde{\varphi}|_{B_p [1,\tau\cdot\sigma]} = \varphi$.
Учитывая предложение 1.8, мы можем считать, что
$\tilde{\varphi}$ --- 
это гомеоморфизм 
$\mathbb{R}^{[1,\tau\cdot\sigma]}$
на
$\mathbb{R}^{[1,\tau\cdot\eta]}$,
являющийся продолжением гомеоморфизма
$\varphi$.
Рассмотрим подпространства 
$M_{\tau\sigma}\subset \mathbb{R}^{[1,\tau\cdot\sigma]}$ 
и
$M_{\tau\eta}\subset \mathbb{R}^{[1,\tau\cdot\eta]}$.
Из леммы 2.5 следует, что 
$\tilde{\varphi}(M_{\tau\sigma}) = M_{\tau\eta}$.

Для каждой точки 
$t\in A_{\tau,\tau\sigma}$
пусть
$\chi_t$ ---
характеристическая функция одноточечного множества
$\{t\}$.
Очевидно, что 
$\{\chi_t\}_{t\in A_{\tau,\tau\sigma}}\subset M_{\tau\sigma}\setminus B_p [1,\tau\cdot\sigma]$
и для любой последовательности попарно различных точек 
$t_n\in A_{\tau,\tau\sigma}$
последовательность 
$\{\chi_{t_n}\}$
поточечно сходится к нулевой функции пространства
$\mathbb{R}^{[1,\tau\cdot\sigma]}$.
Рассмотрим множество функций
$\{\tilde{\varphi}(\chi_t)\}_{t\in A_{\tau,\tau\sigma}}\subset M_{\tau\eta}\setminus B_p [1,\tau\cdot\eta]$.
Каждая из функций 
$\tilde{\varphi}(\chi_t)$
разрывна в некоторой точке множества
$A_{\tau,\tau\eta}\subset [1,\tau\cdot\eta]$,
то есть в некоторой точке вида
$\tau\cdot(\delta +1)$,
где
$\delta < \eta$.
Пусть
$$
B_\delta = \{\tilde{\varphi}(\chi_t)\mid \tilde\varphi(\chi_t) \text{ разрывна в точке } \tau(\delta +1)\}.
$$
Так как 
$|A_{\tau,\tau\eta}|<|\sigma|$
и
$\bigcup_{\delta <\eta}B_\delta = \tilde{\varphi}(\{\chi_t\}_{t\in A_{\tau,\tau\sigma}})$,
то существует такая точка
$\tau(\delta_0 +1)$, $\delta_0 <\eta$,
что 
$|B_{\delta_0}|=|\sigma|>|\eta|>|\delta_0|$.
Поскольку 
$\tilde{\varphi}(\chi_{t_n})\to 0$
для любой последовательности попарно различных точек 
$t_n\in A_{\tau,\tau\sigma}$,
то множество $\{\tilde{\varphi}(\chi_t)\in B_{\delta_0} : \tilde{\varphi}(\chi_t)(\tau(\delta_0+1))\ne 0\}$ не более, чем счётно и, следовательно, множество
$B_{\delta_0}^0 = \{\tilde{\varphi}(\chi_t)\in B_{\delta_0} : \tilde{\varphi}(\chi_t)(\tau(\delta_0+1)) = 0\}$ несчётно.
По лемме 2.3 для каждой функции  
$\tilde{\varphi}(\chi_t)\in B_{\delta_0}^0$ 
существует такой ординал
$\gamma_t$,
что
$\tilde{\varphi}(\chi_t)|_{(\gamma_t,\tau(\delta_0 +1))} = \mathrm{const} = C_t$,
причём
$C_t\ne 0$,
так как все функции из множества
$B_{\delta_0}^0$
разрывны в точке 
$\tau(\delta_0+1)$.
Нетрудно видеть, что существует несчётное число функций из 
$B_{\delta_0}^0$,
для которых 
$|C_t|\geqslant \varepsilon_0$
для некоторого
$\varepsilon_0 >0$.
Рассмотрим последовательность таких попарно различных функций
$\tilde{\varphi}(\chi_{t_n})$,
для которых
$|\tilde{\varphi}(\chi_{t_n})| \equiv C_{t_n}\geqslant\varepsilon_0$ на множестве
$(\gamma_{t_n},\tau (\delta_0 +1))$.
Так как
$\mathrm{cf}(\tau (\delta_0 +1))>\omega$,
то
$\gamma_0 = \sup_{n<\omega}\gamma_{t_n}<\tau (\delta_0 +1)$
и, следовательно,
$|\tilde{\varphi}(\chi_{t_n})(t)|\geqslant\varepsilon_0$
для каждого 
$t\in (\gamma_0, \tau (\delta_0 +1))$.
Это противоречит тому, что последовательность функций
$\{\tilde{\varphi}(\chi_{t_n})\}_{n <\omega}$
поточечно сходится к нулю.
\end{proof}

\textbf{Предложение 2.8. } Пусть $\tau$ --- начальный регулярный ординал, $\tau\geq\omega_2$ и $m,n $ --- различные натуральные числа.  Тогда пространства 
$B_p[1,\tau\cdot m]$ и $B_p[1,\tau\cdot n]$ не гомеоморфны.

Для доказательства этой теоремы нам потребуются некоторые обозначения и вспомогательные утверждения.

Пространство $B_p[1,\tau\cdot n]$ будем отождествлять с пространством $B_p([1,\tau]\times [n])$, где $[n]=\{1,2,\ldots , n\}$. Если $\tau$ --- несчётный регулярный ординал, то для любой функции $x\in B_p[1,\tau]$ найдется такой ординал $\gamma <\tau$, что $x|_{[\gamma,\tau]}$ --- постоянная функция. Отсюда следует, что пространство $B_p[1,\tau]$ линейно гомеоморфно подпространству в $B_p[1,\tau)$, состоящему из всех функций, обращающихся в ноль, начиная с некоторого $\alpha < \tau$. Требуемое линейное гомеоморфное вложение $\varphi\colon B_p[1,\tau]\to B_p[1,\tau)$ можно задать формулой
\begin{align*}
\varphi (x)(1) &= x(\tau);\\
\varphi (x)(n) &= x(n-1)-x(\tau), \text{ если } 2\leqslant n<\omega;\\
\varphi (x)(\alpha) &= x(\alpha)-x(\tau), \text{ если } \omega\leqslant \alpha<\tau.
\end{align*} 
Обозначим это подпространство $B_p^0[1,\tau)$. Очевидно, что и для любого $n\in\mathbb{N}$ пространство $B_p([1,\tau]\times [n])$ линейно гомеоморфно пространству 
$$
B_p^0([1,\tau)\times [n]) = \{x\in B_p([1,\tau)\times [n]) : x|_{[\gamma,\tau)\times [n]}\equiv 0 \text{ для некоторого } \gamma < \tau\}.
$$
Для $\alpha < \tau$ положим
\begin{multline*}
B_p^\alpha([1,\tau)\times [n]) = \\ = \{x\in B_p^0([1,\tau)\times [n]) : x(t,i)=0 \text{ для всех } t\in [\alpha,\tau) \text{ и всех } i\in [n]\}.
\end{multline*}
Если $\alpha$ --- непредельный ординал или $\mathrm{cf}(\alpha) = \omega$, то пространство 
$B_p^\alpha([1,\tau)\times [n])$ можно отождествить с пространством $B_p([1,\alpha)\times [n])$.

Доказательство предложение 2.8 проведем методом от противного. Предположим, что существует гомеоморфизм $T$ между пространствами  $B_p[1,\tau\cdot m]$ и $B_p[1,\tau\cdot n]$. Без ограничения общности можно считать, что $T0=0$.

\textbf{Лемма 2.9.} Пусть $\tau$ --- начальный регулярный ординал, $\tau\geq\omega_2$ и $m,n $ ---  натуральные числа.  Если 
$T\colon B^0_p([1,\tau]\times [m])\to B^0_p([1,\tau]\times [n])$ --- гомеоморфизм, то для любого $\gamma\in (\omega,\tau)$ найдётся $\alpha_\gamma \in(\gamma,\tau)$ такой, что $$
T \left( B_p([1,\alpha_\gamma)\times [m])\right) = B_p([1,\alpha_\gamma)\times [n]).
$$
\begin{proof}
Рассмотрим произвольный ординал $\alpha_1 \in(\gamma,\tau)$. Тогда верно неравенство $d(B_p[1,\alpha_1])\leqslant|[1,\alpha_1]|$, т.е. в $B_p[1,\alpha_1]$ существует всюду плотное подмножество $A$ мощности не больше, чем $|[1,\alpha_1]|$.

Поскольку $\tau$ --- несчётный регулярный ординал, то для любой функции $x\in A$ существует ординал $\beta(x) < \tau$ такой, что $Tx|_{[\beta(x),\tau)\times \{i\}}\equiv 0$ для всех $i=1,2,\ldots ,n$. Тогда 
$\beta_1 = \sup \{\beta(x) : x\in A\}<\tau$ и $Tx|_{[\beta_1,\tau)\times \{i\}}\equiv 0$ для всех 
$x\in B_p([1,\alpha_1]\times [m])$
и для всех $i=1,2,\ldots ,n$.
Можно считать, что $\beta_1 >\alpha_1$.
Аналогично, для ординала $\beta_1 <\tau$
найдётся ординал $\alpha_2$ такой, что $\beta_1 <\alpha_2 <\tau$ и
$T^{-1}(B_p([1,\beta_1]\times [n]))\subset B_p([1,\alpha_2]\times [m])$.
Продолжая этот процесс, получим возрастающую последовательность ординалов 
$$
\gamma <\alpha_1 <\beta_1 <\alpha_2 <\ldots <\alpha_k <\beta_k <\alpha_{k+1}<\ldots <\tau
$$
такую, что
\begin{equation}
T(B_p([1,\alpha_k]\times [m]))\subset B_p([1,\beta_k]\times [n]) \tag{$\ast$}
\end{equation}
и
\begin{equation}
T^{-1}(B_p([1,\beta_k]\times [n]))\subset B_p([1,\alpha_{k+1}]\times [m]). \tag{$\ast\ast$}
\end{equation}
Тогда для ординала $\alpha_\gamma = \sup_{k<\omega}\alpha_k = \sup_{k<\omega}\beta_k$ будет выполнено равенство
$$
T(B_p^{\alpha_\gamma}([1,\tau)\times [m]))=B_p^{\alpha_\gamma}([1,\tau)\times [n]).
$$
Действительно, если $x\in B_p^{\alpha_\gamma}([1,\tau)\times [m])$, то функции $x_k = x|_{[1,\alpha_k]\times [m]}$ являются элементами пространства $B_p([1,\alpha_k]\times [m])$ и сходятся поточечно к $x$. Учитывая включение $(\ast)$, получаем, что 
$Tx_k\in B_p([1,\beta_k]\times [n])\subset B_p^{\alpha_\gamma}([1,\tau)\times [n])$ и, следовательно, $Tx = \lim_{k\to \infty} Tx_k\in B_p^{\alpha_k}([1,\tau)\times [n])$. Обратное включение доказыватся аналогично, с использованием $(\ast\ast)$.
\end{proof}

Так как в лемме 2.9 ординал $\gamma <\tau$ был выбран произвольно, то получаем следующее 

\textbf{Следствие 2.10.} Если
$T\colon B_p^0([1,\tau)\times [m]))\to B_p^0([1,\tau)\times [n])$ --- гомеоморфизм, то множество
$L=\{\alpha\in [1,\tau) : T(B_p^\alpha([1,\tau)\times [m]))=B_p^\alpha([1,\tau)\times [n])\}$
замкнуто и конфинально в промежутке $[1,\tau)$. $\Box$

\textbf{Лемма 2.11.} Пусть $\tau$ --- начальный регулярный ординал, $\tau\geq\omega_2$ и $m,n $ ---  натуральные числа. Если 
$T\colon B_p^0([1,\tau)\times [m])\to B_p^0([1,\tau)\times [n])$ --- гомеоморфизм, то для любого $\gamma <\tau$ найдётся $\alpha_\gamma \in (\gamma,\tau)$ такой, что для функций $x,y\in B_p^0([1,\tau]\times [m])$ условия
$$
x|_{[1,\alpha_\gamma)\times \{i\}} = y|_{[1,\alpha_\gamma)\times \{i\}} \text{ для всех } i=1,2,\ldots ,m
$$
и
$$
Tx|_{[1,\alpha_\gamma)\times \{j\}} = Ty|_{[1,\alpha_\gamma)\times \{j\}} \text{ для всех } j=1,2,\ldots ,n
$$
равносильны. 
\begin{proof}
Рассмотрим произвольный ординал $\alpha_1 <\tau$. Зафиксируем  $t\in [1,\alpha_1]$ и 
натуральное число $k\in [n]$.
Рассмотрим непрерывную функцию $f\colon B_p^0([1,\tau)\times [m])\to\mathbb{R}$, определённую формулой
$f(x) = (Tx)(t,k)$.

Хорошо известно (\cite{Arhangelskii}), что любая непрерывная функция, определенная на всюду плотном подпространстве  в  произведении вещественных прямых, зависит от счетного числа координат, т.е. существует счетное множество $A_{(t,k)}\subset [1,\tau)$ такое, что для функций $x,y\in B_p^0([1,\tau)\times [m])$ из условия 
$x|_{A_{(t,k)}\times [m]} = y|_{A_{(t,k)}\times [m]}$
следует, что
$f(x)=f(y)$.
Отсюда вытекает, что найдется такой ординал $\beta_1\in (\alpha_1,\tau)$, для которого из того, что
$x|_{[1,\beta_1]\times [m]} = y|_{[1,\beta_1]\times [m]}$
следует равенство
$Tx|_{[1,\alpha_1]\times [n]} = Ty|_{[1,\alpha_1]\times [n]}$.

Аналогично, для ординала $\beta_1$ найдётся ординал $\alpha_2\in (\beta_1,\tau)$ такой, что из условия 
$Tx|_{[1,\alpha_2]\times [n]} = Ty|_{[1,\alpha_2]\times [n]}$ 
следует равенство 
$x|_{[1,\beta_1]\times [m]} = y|_{[1,\beta_1]\times [m]}$.

Продолжая этот процесс, получим последовательность
$$
\gamma <\alpha_1 <\beta_1 <\alpha_2 <\ldots <\alpha_k <\beta_k <\alpha_{k+1}<\ldots <\tau,
$$
для которой условие
$x|_{[1,\beta_k]\times [m]} = y|_{[1,\beta_k]\times [m]}$
влечёт 
$Tx|_{[1,\alpha_k]\times [n]} = Ty|_{[1,\alpha_k]\times [n]}$
и условие
$Tx|_{[1,\alpha_{k+1}]\times [n]} = Ty|_{[1,\alpha_{k+1}]\times [n]}$ 
влечёт  
$x|_{[1,\beta_k]\times [m]} = y|_{[1,\beta_k]\times [m]}$.

Тогда для ординала $\alpha_\gamma = \sup_{k<\omega}\alpha_k = \sup_{k<\omega}\beta_k$ будет верно утверждение леммы.
\end{proof}

\textbf{Следствие 2.12.} Если
$T\colon B_p^0([1,\tau)\times [m])\to B_p^0([1,\tau)\times [n])$ --- гомеоморфизм, то множество
$$
M=\{\alpha\in [1,\tau) : x|_{[1,\alpha)\times [m]} = y|_{[1,\alpha)\times [m]} \Longleftrightarrow Tx|_{[1,\alpha)\times [n]} = Ty|_{[1,\alpha)\times [n]}\}
$$ 
замкнуто и конфинально в промежутке $[1,\tau)$. $\Box$

\textbf{Доказательство предложения 2.8.} Предположим, что существует гомеоморфизм 
$T\colon B_p^0([1,\tau)\times [m])\to B_p^0([1,\tau)\times [n])$. Поскольку $\tau\geqslant \omega_2$ и множества $L$ и $M$ конфинальны и замкнуты в $[1,\tau)$, то найдется такой ординал $\gamma_0\in L\cap M$, что $\mathrm{cf}(\gamma_0)=\omega_1$. 

Рассмотрим в $B_p^0([1,\tau)\times [m])$ и $B_p^0([1,\tau)\times [n])$ замкнутые подпространства
$$
E=\{x\in B_p^0([1,\tau)\times [m]) : x(\gamma_0,i)=0, \ i=1,2,\ldots , m \}
$$ 
и
$$
H=\{x\in B_p^0([1,\tau)\times [n]) : x(\gamma_0,i)=0, \ i=1,2,\ldots , n \}
$$
соответственно. 
Для $x\in E$ пусть 
$$
x'(t,i)=\left\lbrace 
\begin{array}{ll}
x(t,i), & \text{ если } t\leqslant\gamma_0;\\
0, & \text{ если } t >\gamma_0.
\end{array}\right. 
$$
Так как 
$\gamma_0\in L$, то $T(x')\in B_p^0([1,\gamma_0)\times [n])$
и, значит,
$T(x')(\gamma_0,i)=0$ для $i=1,2,\ldots , n$.
С другой стороны, так как 
$x'|_{[1,\gamma_0)\times [m]} = x|_{[1,\gamma_0)\times [m]}$
и 
$\gamma_0\in M$,
то
$T(x')|_{[1,\gamma_0)\times [n]} = T(x)|_{[1,\gamma_0)\times [n]}$
и, в силу непрерывности функций 
$T(x)$ и $T(x')$
в точке $\gamma_0$
получаем, что
$T(x')(\gamma_0,i) = T(x)(\gamma_0,i)$ для всех $i=1,2,\ldots ,n$.
Следовательно, 
$T(x)(\gamma_0,i)=0$ для всех $i=1,2,\ldots ,n$, 
откуда $T(E)\subset H$.
Аналогично получаем, что $T^{-1}(H)\subset E$, то есть $T(E)=H$. Но это невозможно, так как подпространства $E$ и $H$ имеют различные коразмерности в $B_p^0([1,\tau)\times [m])$ и $B_p^0([1,\tau)\times [n])$ соответственно (см. лемму 4 в \cite{Gulko3}). 

\textbf{Предложение 2.13.} Пусть $\tau\geqslant\omega_2$ --- регулярный ординал и $n <\omega$. Тогда пространства $B_p[1,\tau\cdot\omega]$ и $B_p[1,\tau\cdot n]$ не гомеоморфны.
\begin{proof}
Это утверждение становится очевидным, если заметить, что пространство $B_p[1,\tau\cdot\omega]$ гомеоморфно своему квадрату, а пространство $B_p[1,\tau\cdot n]$ своему квадрату не гомеоморфно по предложению 2.8.
\end{proof}

Предложения 2.8 и 2.13 окончательно завершают доказательство теоремы 2.2.

Авторы выражают глубокую признательность рецензенту за внимательное отношение к работе и полезные замечания.

\medskip

Гензе Леонид Владимирович

Томский государственный университет

Механико-математический факультет

634050 г.~Томск, пр-т Ленина, 36

genze@math.tsu.ru 

\bigskip

Гулько Сергей Порфирьевич

Томский государственный университет

Механико-математический факультет

634050 г.~Томск, пр-т Ленина, 36

gulko@math.tsu.ru

\bigskip

Хмылева Татьяна Евгеньевна

Томский государственный университет

Механико-математический факультет

634050 г.~Томск, пр-т Ленина, 36

tex2150@yandex.ru


\medskip

\begin{thebibliography}{100}

\bibitem{BesPel}
{\it Bessaga~C., Pelczynski~A.} On isomorphic classification
of spaces of continuos functions. Studia Math. 1960. V.~19. P.~53--62.

\bibitem{Semadeni}
{\it Semadeni~Z.} Banach spaces non-isomorphic to their Cartesian product.
Bull. Acad. Pol. Sci. Ser. Math. Stron. et Phys. 1960. V.~8. P.~81--84.

\bibitem{GulOs}
{\it Гулько~C.\,П., Оськин~А.\,В.} Изоморфная классификация пространств непрерывных функций на вполне упорядочен ных бикомпактах .Функциональный анализ и его приложения.1975. Т.~9. №~1. С.~61-62.

\bibitem{Kisl}
{\it Кисляков~С.\,В.} Изоморфная классификация пространств непрерывных функций на ординалах. Сиб. мат. журнал.
1975. Т.~16. С.~293-300.

\bibitem{Gulko}
{\it Гулько~С.\,П.} Свободные топологические группы и
пространства непрерывных функций на ординалах. Вестник томского государственного унив-та. 2003. №~280.

\bibitem{BaarsDeGroot}
{\it Baars J., de Groot J.} On topological and linear equivalence of certain function spaces, CWI Tract 86, Stichting Mathematisch Centrum, Centrum voor Wiskunde en Informatica, Amsterdam, 1992.

\bibitem{GGK}
\textit{Гензе Л.В., Гулько С.П., Хмылева Т.Е.} Классификация пространств бэровских функций на отрезках ординалов // Труды Института математики и механики УрО РАН. 2010. Т. 16. № 3. С. 61--66.


\bibitem{Engel}  
\textit{Энгелькинг Р.} Общая топология. М.: Мир, 1986.

\bibitem{Arch}
{\it Архангельский~А.В.} О линейных гомеоморфизмах пространств функций. ДАН СССР. 1982. Т.~264, №~6. С.~1289 – 1292.

\bibitem{Tkachuk}
{\it Tkachuk~V.V.} A $C_p$--- Theory Problem Book. Topological and Function Spaces. Springer. 2011. 486pp.

\bibitem{Arhangelskii}
{\it Архангельский~А.В.} Топологические пространcтва функций. М.:
Издательство МГУ. 1989.

\bibitem{Gulko3}
{\it Гулько~С.П.} Пространства непрерывных функций на ординалах и ультрафильтрах // Матем. заметки, 47:4 (1990), 26–34.

\end{thebibliography}
\end{document}